\documentclass[11pt]{amsart}
\usepackage{geometry}                % See geometry.pdf to learn the layout options. There are lots.
\usepackage{graphicx}
\usepackage{amssymb}
\usepackage{epstopdf}
\newcommand{\C}{\mathbb{C}}
\newcommand{\R}{\mathbb{R}}
\newcommand{\Z}{\mathbb{Z}}
\newcommand{\eps}{\varepsilon}
\newcommand{\tr}{\mbox{tr}}

\title{Towards a Spectral Proof of Riemann's Hypothesis}
\author{R.S.MacKay}
\address{Mathematics Institute and Centre for Complexity Science,
University of Warwick, Coventry CV4 7AL, UK}
\date{\today}                                           % Activate to display a given date or no date

\begin{document}
\begin{abstract}
The paper presents evidence that Riemann's $\xi$ function evaluated at $2\sqrt{E}$ could be the characteristic function $P(E)$ for the magnetic Laplacian minus $\frac{85}{16}$ on a surface of curvature $-1$ with magnetic field $\frac94$, a cusp of width 1, a Dirichlet condition at a point, and other conditions not yet determined.
\end{abstract}
\maketitle

\section{Introduction}
Riemann's hypothesis (RH) is that all the zeroes of his $\xi$-function are real \cite{R}.  He expressed $\xi$ as 
$$\xi(\omega) = 4 \int_1^\infty \frac{d(x^{3/2}\psi'(x))}{dx} x^{-1/4} \cos(\frac{\omega}{2} \log{x})\ dx,$$
where $\psi(x) = \sum_{n=1}^\infty e^{-n^2\pi x}$, but a change of variables makes it the Fourier transform of a function that Polya called $\Phi$ \cite{Po1}:
\begin{equation}
\xi(\omega) = \int_{-\infty}^\infty e^{i\omega t} \Phi(t)\ dt = 2 \int_0^\infty \cos \omega t \ \Phi(t)\ dt,
\label{eq:xiPhi}
\end{equation}
with
\begin{equation}
\Phi(t) = \sum_{n\ge 1} (4 \pi^2 n^4 e^{9t/2} - 6\pi n^2 e^{5t/2})e^{-\pi n^2 e^{2t}}.
\end{equation}
Note that:
\begin{enumerate}
\item I denote the argument of $\xi$ by $\omega$ because I wish to think of it as a frequency;
\item Confusingly, many authors (supported by \cite{Ed}) write $\xi(s)$ for $\xi(i(\frac12-s))$;
\item $\Phi$ is even, because $\Phi(t)$ can be written as $\frac12 (\partial_t^2-\frac14)(e^{t/2}\theta(e^{2t}))$ where $\theta(v) = \sum_{n\in \Z} e^{-\pi n^2 v} = 1 + 2 \psi(v)$ is Jacobi's $\theta$ function, which satisfies $v^{1/4} \theta(v) = v^{-1/4} \theta(1/v)$ (beware, some people write $\theta(u)$ for $\theta(-iu)$).
\end{enumerate}
It follows that $\xi$ is entire and even.  Furthermore it is of order 1 (for given modulus of $\omega$, $|\xi(\omega)|$ is maximised on the imaginary axis). 

The $\xi$ function is related to the $\zeta$ function, defined for $\Re s > 1$ by
\begin{equation}
\zeta(s) = \sum_{n\ge 1} n^{-s} = \prod_{p \mbox{ prime}} (1-p^{-s})^{-1},
\end{equation}
via \cite{R}
\begin{equation}
\xi(\omega) = -\frac12(\omega^2+\frac14) \Gamma(\frac14+i\frac{\omega}{2}) \pi^{-\frac14-i\frac{\omega}{2}} \zeta(\frac12+ i\omega).
\label{eq:xizeta}
\end{equation}

The prime number theorem \cite{Had,VP} states that that the number $\pi(x)$ of primes less than $x$ satisfies 
$$\pi(x) = Li(x) + R(x)$$ 
with $Li(x) = \int^x \frac{dy}{\log y} \sim \frac{x}{\log x}$ as $x \to +\infty$ (conventions differ for the lower limit of integration but they make no difference to these results) and a remainder $R(x)$ whose size relative to $Li(x)$ goes to zero as $x \to \infty$.
If RH is true then the remainder term in the prime number theorem can be improved from the best currently known \cite{Wa}
\begin{equation}
R(x) = O\left(x \exp{-\frac{A (\log x)^{3/5}}{(\log\log x)^{1/5}}}\right)
\end{equation}
to \cite{Koch}
\begin{equation}
R(x) = O(\sqrt{x} \log x).
\end{equation}

A strategy to prove RH is to show that the zeroes of $\xi$ correspond to the eigenvalues of some Hermitian operator $H$, because the spectrum of any Hermitian operator is real.  This strategy has been attributed to Hilbert and Polya, supported by the letter \cite{Po2}\footnote{though it is strange that Polya made no mention of the spectral strategy in his paper \cite{Po1} on the subject of RH, published nine years after the claimed suggestion.  Furthermore, he did not even use the spectral strategy there to prove that his function $2\mathfrak{G}(i\frac{\omega}{2};\pi)$, obtained by replacing $\Phi(t)$ by $e^{-2\pi \cosh 2t}$, has all its zeroes real, despite its being recognisable \cite{Bi} as the modified Bessel function (or MacDonald function) $K_{i\omega/2}(2\pi)$, which is zero iff $\omega^2/4$ is an eigenvalue of the modified Bessel operator $H = -z^2\partial_z^2 - z\partial_z + z^2$ on functions on $[2\pi,\infty)$ with Dirichlet boundary conditions ($H$ is Hermitian with respect to the scalar product $\langle \phi, \psi \rangle = \int_{2\pi}^\infty \bar{\phi}(z) \psi(z) dz/z$).}.

Although the simplest interpretation of ``correspond'' is that the zeroes of $\xi$ are the eigenvalues of $H$, evenness of $\xi$ makes it natural to seek a positive Hermitian operator whose eigenvalues are the squares of the Riemann zeroes, which would also prove RH.  Think of classical acoustic or membrane problems where the eigenvalues are the squares of the frequencies of modes of vibration.
Actually, it is not necessary to require $H$ to be positive as it is already established that $\xi$ has no purely imaginary zeroes.

This is the approach that I shall explore\footnote{though I admit that it might be fruitful instead to seek a Hermitian operator of Dirac form, e.g.~$(z\partial_z + g(z))\psi_1 = \omega \psi_2$, $(-z\partial_z + g(z))\psi_2 = \omega \psi_1$, whose eigenvalues $\omega$ would be the Riemann zeroes in $\pm$ pairs.  Another alternative is to seek a linear positive-definite Hamiltonian system whose eigenvalues are $i$ times the Riemann zeroes, because the eigenvalues of such a system come in $\pm$ pairs and are pure imaginary (Dirichlet's criterion, e.g.~\cite{Mac}), a suggestion I made to Michael Berry in 1997.}.
Thus, let\footnote{the factor 2 is just to fit with the Bessel operator  mentioned in the first footnote.} 
\begin{equation}
\Xi(E) = \xi(2\sqrt{E})
\end{equation}
and seek a Hermitian operator whose eigenvalues are the zeroes of $\Xi$.  Note that $\Xi$ is entire of order $\frac12$.

It is not just the zeroes of $\Xi$ that one may hope to be the eigenvalues of a Hermitian operator, but the whole function $\Xi$ could be the characteristic function of the operator.  By ``characteristic function'' I mean an analytic function $P$ whose zeroes are the eigenvalues.  Constructions of such functions go under the names of spectral determinant and functional determinant and can involve zeta-regularisation or functional integration.  Definitions can differ by multiplication by a nowhere-zero entire function, but such a function is the exponential of an entire function so if not a constant it has order at least 1.  The order of the product of entire functions of different orders is the larger of the two orders \cite{Lev}.  Thus if there is an order $\frac12$ characteristic function then it is unique up to multiplication by a non-zero constant.  

If the resolvent operator $(E-H)^{-1}$ is of trace-class, my preferred definition for the characteristic function is
\begin{equation}
P(E) = C \exp \int_{E_0}^E R(E')\ dE'
\end{equation}
where 
\begin{equation}
R(E) = \tr (E-H)^{-1},
\end{equation}
$E_0$ is an arbitrary point not in the spectrum, and $C$ is an arbitrary non-zero constant.  The function $P$ is well defined and entire because although $R$ has a pole at each eigenvalue, the residue is an integer (the algebraic multiplicity of the eigenvalue), so its contribution to the exponential is a factor of 1, and the value of $P$ at an eigenvalue can be filled in analytically by $0$.

Thus I aim to find a Hermitian operator $H$ whose characteristic function is $\Xi$, or equivalently, for which $R(E) = \partial_E \log \Xi(E)$.  The goal is not achieved here, but I report on progress.

\section{1D operators}
The simplest class of Hermitian operators to try are the one-dimensional Schr\"odinger operators\footnote{one might propose to enlarge the scope to Sturm-Liouville operators but they can all be transformed to Schr\"odinger form.  To spell this out, the general Sturm-Liouville operator is $Ly = \frac{1}{w}(-(ay')'+dy)$ on $I=\R$ or a subinterval, where $w$ and $a$ are positive functions and $d$ goes to infinity at open ends of $I$.  It is Hermitian with respect to the scalar product $\langle y_1, y_2 \rangle = \int_I \bar{y}_1(z) y_2(z) w(z)\ dz$.  Let $x = \int^z \sqrt{w(z)/a(z)}\ dz$, $Q = \log aw$ and write $y = e^{-Q/4} \psi$.  Then $e^{Q/4}Le^{-Q/4} $ is of Schr\"odinger form with $V(x) = \frac14 Q''(z) + \frac{1}{16} Q'^2(z) + d(z)e^{-Q(z)/2}$.}.  
These are operators of the form
\begin{equation}
H = -\partial_x^2 + V(x)
\end{equation}
on functions $\psi$ from $I=\R$ or a subinterval to $\R$ and we take Dirichlet boundary conditions $\psi(x) \to 0$ at the ends.  The ``potential'' $V$ is assumed to go to $+\infty$ at open ends of $I$, to be locally integrable and to have $\int_{V\le 0} \sqrt{-V(x)} \ dx < \infty$.  Then there is a subspace of $L^2(I)$ on which $H$ is Hermitian with respect to the standard scalar product, and its spectrum is real, discrete and bounded below. 
%[REF]

We suspect in advance that this class will fail to contain a match to the Riemann zeroes, because:
\begin{enumerate}
\item the spacing of eigenvalues of a (smooth enough) 1D Schr\"odinger operator is too regular compared to numerics for the Riemann zeroes;
\item the statistics of the spacings between Riemann zeroes suggest a complex Hermitian operator, whereas 1D Schr\"odinger operators are real.
\end{enumerate}
Nevertheless, I believe it is an informative starting point.

First, let us use the known asymptotics of the number $N(\Omega)$ of Riemann zeroes with real part in $[0,\Omega]$ to determine the asymptotics of the potential $V$.  Riemann \cite{R} claimed by integrating $d\log \xi(t)$ along a contour around the rectangle $|\Im \omega| \le \frac12, \Re \omega \in [0,\Omega]$, that
\begin{equation}
N(\Omega) = \frac{\Omega}{2\pi} \log \frac{\Omega}{2\pi e} + O(\log \Omega).
\end{equation}
It is not clear whether he had a proof but one was eventually provided by von Mangoldt \cite{Ma} (a simpler proof was given later in \cite{Ba}).  Because $|\Im \omega| \le \frac12$ for the zeroes of $\xi$, we obtain the number of zeroes of $\Xi$ with real part less than $E$ 
\begin{equation}
N_\Xi(E) = \frac{\sqrt{E}}{\pi} \log \frac{\sqrt{E}}{\pi e} + O(\log E)
\end{equation}
(it is sandwiched between $N(2\sqrt{E})$ and $N(2\sqrt{E+\frac{1}{16}})$).

From semiclassical quantisation theory, 
%[REF], 
the number $N(E)$ of eigenvalues of a Schr\"odinger operator less than $E$ is asymptotically $W(E)/2\pi$, where
\begin{equation}
W(E) = 2 \int_{V(x)\le E} \sqrt{E-V(x)}\ dx.
\end{equation}
The error term is $O(E^{-\frac12})$.
Define the {\em width} of the potential at height $v$ to be
\begin{equation}
w(v) = \int_{V(x)\le v} dx.
\end{equation}
Then by breaking $I$ into intervals of monotonicity of $V$,
\begin{equation}
W(E) = 2\int_{-\infty}^E \sqrt{E-v}\ dw(v),
\end{equation}
a Riemann-Stieltjes integral\footnote{one can think of $dw(v) = w'(v) dv$ with $w'$ as a distributional derivative; note that if $V$ has infimum $V_{min}$ then one could start the integral at $V_{min}$ but $w(v)=0$ for $v < V_{min}$ so we start it at $-\infty$,$w'(v) \to \infty$ from below or above as $v$ approaches local maxima or minima respectively of $V$, and if $V$ has a plateau then $w$ has a jump.}.  
%So
%\begin{equation}
%N(E) \sim \frac{1}{\pi} \int_{-\infty}^E \sqrt{E-v}\ dw(v).
%\end{equation}

We want to determine $w(v)$, given the asymptotics of $W(E)$.  This is achieved by Abel inversion (e.g.~\cite{Ke}): 
\begin{equation}
w(v) = \frac{1}{\pi} \int_{-\infty}^v \frac{dW(E)}{\sqrt{v-E}}.
\end{equation}

Take $W(E) = 2\sqrt{E} \log \frac{\sqrt{E}}{\pi e}$ for $E > 0$, and $0$ for $E\le 0$ (it is continuous at 0), being a leading  approximation to $2\pi N_\Xi(E)$.  Then $W'(E) = \frac{1}{\sqrt{E}} \log \frac{\sqrt{E}}{\pi}$ for $E > 0$, and 0 for $E < 0$, with no delta-function at 0.  So
\begin{equation}
w(v) = \frac{1}{\pi} \int_0^v \frac{\log\frac{\sqrt{E}}{\pi}}{\sqrt{E}\sqrt{v-E}}\ dE.
\end{equation}
Write $E=v \sin^2\theta$ and change variable of integration to $\theta$:
\begin{equation}
w(v) = \frac{1}{\pi} \log \frac{\sqrt{v}}{\pi} \int_0^{\pi/2} 2\ d\theta + \frac{1}{\pi} \int_0^{\pi/2} 2\log\sin \theta \ d\theta.
\end{equation}
The first integral evaluates to $\pi$ and the second\footnote{One way to derive this integral is to write $\sin \theta = 2 \sin \frac{\theta}{2} \cos \frac{\theta}{2}$ and use $\cos \frac{\theta}{2} = \sin\frac{\pi-\theta}{2}$ to see that $L : = \int_0^{\pi/2} 2 \log \sin \theta\ d\theta = \pi \log 2 + 2L$.} to $-\pi \log 2$.
Thus for this approximation to $W$,
\begin{equation}
w(v) = \log \frac{\sqrt{v}}{2\pi}.
\end{equation}

If we allow a correction $O(\log E)$ to $W(E)$, Abel transforming it produces $O(\frac{1}{\sqrt{v}}\log v)$ (by a similar calculation), and the Abel transform is linear so
\begin{equation}
w(v) = \log \frac{\sqrt{v}}{2\pi} + O(\frac{1}{\sqrt{v}}\log v).
\end{equation}
Inverting, we obtain
\begin{equation}
v(w) = 4\pi^2 e^{2w} + O(w e^w)
\label{eq:vasymp}
\end{equation}
as $w \to \infty$.\footnote{Actually, I suspect one should be allowed to match the ``smoothed'' number of Riemann zeroes
$$\bar{N}(\Omega) = \frac{\Omega}{2\pi} \log \frac{\Omega}{2\pi e} + \frac78 + O(\Omega^{-1})$$
with $W(E)$ ($E=\Omega^2/4$).  In this case, the correction to $W$ would be $O(E^{-\frac12})$, that to $w$ would be $O(v^{-\frac12})$ and that to $v$ would be $O(e^w)$.}

There are many potentials $V$ with the asymptotics (\ref{eq:vasymp}).  The simplest is to take one-sided potentials with $V(x) \sim 4\pi^2 e^{2x}$ for $x > 0$ and a hard wall at $x=0$ (i.e.~Dirichlet boundary condition $\psi(0)=0$).  The special case $V(x) = 4\pi^2 e^{2x}$ is solvable in terms of modified Bessel functions by the change of variable $z = 2\pi e^x$.  For energy $E$, the solutions decaying at infinity are the multiples of $K_{\sqrt{-E}}(2\pi e^x)$.\footnote{$K_\nu$ is even in its order $\nu$ thus the square root induces no singularity.}  Thus the condition for $E$ to be an eigenvalue is that it be a zero of $K_{\sqrt{-E}}(2\pi)$, so we define its characteristic function $P_K(E) = K_{\sqrt{-E}} (2\pi)$.  This is Polya's $2\mathfrak{G}(i\frac{\omega}{2};\pi)$ \cite{Po1}, reached by a different route.  Note that it has order $\frac12$ as required.

\begin{figure}[htbp] %  figure placement: here, top, bottom, or page
   \centering
   \includegraphics[width=5in]{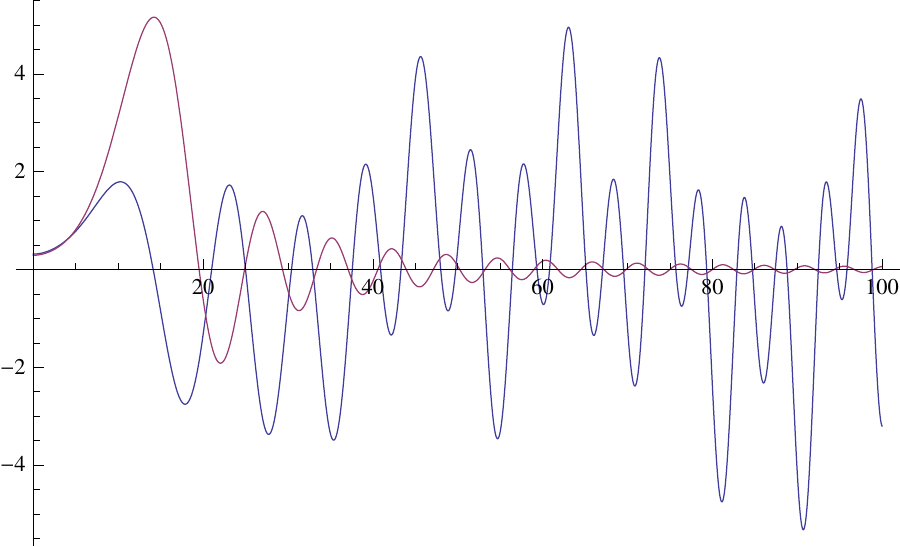} 
   \caption{$S(\omega)\xi(\omega)$ (in blue) and $500 S(\omega) K_{i\frac{\omega}{2}}(2\pi)$ (in red) for $\omega \in [0,100]$.}
   \label{fig:xiK}
\end{figure}

To compare the graphs of the two functions $P_K$ and $\Xi$, it is necessary to scale them, because they are huge for large $E<0$ and tiny for large $E>0$.  It is conventional to factorise $\xi(\omega) = -f(\omega)Z(\omega)$ into real analytic even functions:
\begin{eqnarray}
f(\omega) &=& \frac12 \pi^{-1/4} (\omega^2 + \frac14) \sqrt{\Gamma(\frac14+i\frac{\omega}{2})\Gamma(\frac14-i\frac{\omega}{2})} \\
Z(\omega) &=& \sqrt{\frac{\Gamma(\frac14+i\frac{\omega}{2})}{\Gamma(\frac14- i \frac{\omega}{2})}}\pi^{-i\omega/2}\zeta(\frac12+i\omega) ,
\end{eqnarray}
of which $f$ is positive for real $\omega$, and then to plot $Z$, which is $-\xi/f$.
I prefer to define
\begin{equation}
S(\omega) = \frac{2^{3/2} \cosh \frac{\pi \omega}{4}}{\pi^{1/4} (\omega^2 + 4)^{7/8}},
\end{equation}
which is asymptotically the same as $1/f(\omega)$, and plot $S\xi$, because it is smoother than $1/f$ at $\omega=0$.\footnote{$Z$ has a blip there coming from its poles at $\omega= \pm i/2$, whereas the closest poles of $S\xi$ are at $\omega = \pm 2i$.} Thus Figure~\ref{fig:xiK} shows $S(\omega)\xi(\omega)$ and $500 S(\omega) K_{i\frac{\omega}{2}}(2\pi)$ for a range of real $\omega$.\footnote{500 is a convenient numerical factor to bring the two functions to roughly the same value at 0; one could work out the exact multiple to make them agree at 0.}  The fit is not good.

Examination of $P_K$ and $\Xi$ for negative $E$ reveals a poor fit there too.  Writing $E = -\nu^2/4$ for negative $E$ and using (\ref{eq:xizeta}), 
\begin{equation}
\xi(-i\nu) \sim (\frac{\nu}{2\pi e})^{\frac{\nu}{2}} \nu^{7/4} (\frac{\pi}{2})^\frac14,
\end{equation}
whereas from the integral representation $K_\nu(z) = \int_0^\infty \cosh \nu t \ e^{-z\cosh t} dt$, one obtains
\begin{equation}
K_{\nu/2}(2\pi) \sim \sqrt{\frac{\pi}{\nu}} (\frac{\nu}{2\pi e})^{\frac{\nu}{2}} ,
\end{equation}
which is too small by a factor  proportional to $\nu^{9/4}$.

The situation cannot be saved by considering two-sided potentials with the same width function, because as I shall show shortly, the negative $E$ asymptotics is determined by the width function and the number of sides of the potential, and the two-sided case does not match $\Xi$ either.  It is worth commenting, however, that the two-sided case $V(x) = 4\pi^2 e^{4|x|}$ is solvable again by modified Bessel functions and $V(x) = 8\pi^2 \cosh 4x$ is solvable by modified Mathieu functions.  There are many asymmetric potentials with the same asymptotic width too. 
One example is the Tzitzecka potential $V(x) = 4\pi^2 2^{-2/3} (2 e^{3x} + e^{-6x})$.
%I do not know one solvable in terms of standard functions, though that is not necessary for their analysis, just a convenient short cut.

It is also worth commenting that Polya \cite{Po1} proved by an ingenious argument that the Fourier transform of $2\pi^2 \cosh \frac{9t}{2} e^{-2\pi \cosh 2t}$, which is a good even approximation to $\Phi(t)$ for large $t$, has only real zeroes.  He called it the fake $\xi$-function $\xi^*$.  It has the right asymptotics for negative $E$.  In an afterword, he wrote that one could do the same for the better approximation $(2\pi^2 \cosh \frac{9t}{2} - 3 \pi \cosh \frac{5t}{2}) e^{-2\pi \cosh 2t}$.

But how about finding a potential $V$ for which the Schr\"odinger operator has characteristic function with the right negative $E$ asymptotics?  

For 1D Schr\"odinger operators, the characteristic function has a simple expression in terms of shooting, which we have already seen for the case of the modified Bessel operator.  This is said by various authors to have been explained by Gelfand and Yaglom \cite{GY}.
%\footnote{though I can't find where in the paper it is done.}. 
First consider the case of a potential with $V(x) \to +\infty$ as $x \to +\infty$ and a hard wall, without loss of generality at $x=0$.  Then for any $E \in \C$ there is a solution $\psi_E^+$ going to 0 at infinity, unique up to scalar multiple.  Choose a reference value of $E$, for example 0, and a solution $\psi_0^+$ going to 0 at infinity.  Then the ratio $\psi_E^+/\psi_0^+$ goes to a limit as $x \to \infty$, so we can normalise $\psi_E^+$ to be asymptotic to $\psi_0^+$.  The characteristic function for the operator is then $P(E) = \psi_E^+(0)$, because this is zero iff $E$ is an eigenvalue.  The same applies if $V(x) \to \infty$ as $x$ approaches a finite value $a>0$ from the left.

Next consider two-sided potentials.  I describe this in the case that $V(x) \to +\infty$ as $x \to \pm \infty$ but again one can allow finite limit points.  Then for reference energy 0, choose solutions $\psi_0^\pm$ going to zero to the right and left respectively.  For general $E$, let $\psi_E^\pm$ be the solutions asymptotic to $\psi_0^\pm$ to the right and left respectively.  Then $P(E)$ is the Wronskian of the two solutions, evaluated anywhere as it is constant, but for definiteness let us say at $x=0$: $P(E) = \psi_E^+(0) \partial_x \psi_E^-(0) - \psi_E^-(0) \partial_x \psi_E^+(0)$.  This is because for a second order ordinary differential equation, $\psi^+$ is a multiple of $\psi^-$ (making an eigenfunction) iff the Wronskian is zero.

Almost as for the asymptotics of the number of eigenvalues, the asymptotics of the characteristic function for negative $E$ depends only on the width function of the potential and whether it is one- or two-sided.  Here is the analysis.

Using the Liouville-Green-Wentzell-Kramers-Brillouin (LGWKB) method,
%[REF], 
assuming $V$ is bounded below and $V' \ll (V-E)^\frac32$,
%\footnote{no need for Jeffreys' contribution here, because that showed how to take care of turning points, which do not occur for large negative $E$}, 
\begin{equation}
\psi_E^+(0) \sim C(V(0)-E)^{-\frac14} \exp{\int_0^\infty (\sqrt{V(x)-E}-\sqrt{V(x)})\ dx},
\end{equation}
where the constant $C$ is related to the normalisation of $\psi_0^+$.  To remove its effect it is convenient to take the logarithmic derivative.  Differentiation of the above asymptotics with respect to $E$ can be justified, so in the one-sided case we obtain
\begin{equation}
R(E) = \partial_E \log P(E) \sim -\frac{1}{4E} - \frac12 \int_0^\infty \frac{dx}{\sqrt{V(x)-E}}.
\label{eq:R1}
\end{equation}

In the two-sided case, we need the $x$-derivative of $\psi_E^+$:
\begin{equation}
\partial_x \psi_E^+(0) \sim - (V(0)-E)^\frac14 \exp{\int_0^\infty (\sqrt{V(x)-E}-\sqrt{V(x)})\ dx} .
\end{equation}
Doing the same for $\psi_E^-$, we obtain
\begin{equation}
P(E) \sim \exp{\int_{-\infty}^\infty (\sqrt{V(x)-E}-\sqrt{V(x)})\ dx} .
\end{equation}
Thus
\begin{equation}
R(E) \sim -\frac12 \int_{-\infty}^\infty \frac{dx}{\sqrt{V(x)-E}} .
\label{eq:R2}
\end{equation}

I denote 
\begin{equation}
T(E) = \frac12 \int_{-\infty}^\infty \frac{dx}{\sqrt{V(x)-E}},
\end{equation}
because it is the imaginary time for a classical particle of mass $\frac12$ to pass under the potential at energy $E$.  It depends only on the width function for the potential:
\begin{equation}
T(E) = \frac12 \int_{-\infty}^\infty \frac{dw(v)}{\sqrt{v-E}} = \frac12 \int_0^\infty \frac{dw}{\sqrt{v(w)-E}},
\end{equation}
where $v(w)$ is the inverse function to $w(v)$.

To compare with (\ref{eq:R1}) or (\ref{eq:R2}), from (\ref{eq:xizeta}) and using the asymptotic expansion of the digamma function 
\begin{equation}
\Psi(z) = (\log \Gamma)'(z) \sim \log z - \frac{1}{2z} + O(z^{-2})
\end{equation}
as $z \to  +\infty$, one can obtain
\begin{equation}
\partial_E \log \Xi(E) = -\frac{1}{2\sqrt{-E}} \log \frac{\sqrt{-E}}{\pi} + \frac{7}{8E} + O((-E)^{-\frac32}).
\end{equation}
for large negative $E$.

We are faced with another inverse problem: given the function $T$, deduce the function $w$.  We want $w$ so that 
\begin{equation}
T(E) \sim \frac{1}{2\sqrt{-E}} \log \frac{\sqrt{-E}}{\pi} - \frac{\kappa}{2E},
\label{eq:T}
\end{equation}
where $\kappa=7/4$ for a two-sided potential, $9/4$ for a one-sided potential.

I did not manage to make an inversion formula (in contrast to the Abel case).  But taking
\begin{equation}
v(w) = 4\pi^2 e^{2w} - \beta e^w + \gamma
\end{equation}
(with $\beta \le 8\pi^2$ to ensure monotonicity, which the inverse of any width function must satisfy) produces
\begin{equation}
T(E) = \frac{1}{2\sqrt{\gamma-E}} \log \frac{2\sqrt{\gamma-E}\sqrt{\gamma-E+4\pi^2 - \beta} + 2(\gamma-E)-\beta} {4\pi \sqrt{\gamma-E}-\beta}.
\end{equation}
This fits with (\ref{eq:T}) to order $E^{-1}$ iff $$\beta = 4 \pi \kappa.$$
The effects of $\gamma$ are below the accuracy of the above asymptotics for the shooting function, so to determine its value would require further accuracy.  

For the one-sided case, this gives the Morse potential 
\begin{equation}
V(x) = 4\pi^2 e^{2x} - 4\pi\kappa e^x + \gamma
\label{eq:Morse}
\end{equation}
with hard wall at $x=0$ and $\kappa = 9/4$.  The Morse potential is usually studied on the whole line rather than the half-line, as an approximation to chemical bonds.  Its possible relevance to RH was noticed also by \cite{Lag}, but without identification of the required value of $\beta$.\footnote{he highlights $\beta = 2\pi$ as a special case ($k=-\frac12$ in his notation), but he does not even determine the leading coefficient $4\pi^2$, though it is implicit in his results.}

The Schr\"odinger equation 
\begin{equation}
(-\partial_x^2 + 4\pi e^{2x} - 4\pi \kappa e^x+\gamma) \psi = E\psi
\end{equation}
transforms to Whittaker's equation 
\begin{equation}
(-z^2\partial_z^2 + \frac{z^2-1}{4} - \kappa z)\phi = -\mu^2 \phi
\end{equation}
under $z = 4\pi e^x$ (the extra factor of 2 is a quirk of Whittaker's conventions), $\psi = z^{-1/2} \phi$ and $E = \gamma-\mu^2$.   
Its solutions decaying as $z \to +\infty$ are the multiples of the Whittaker function $W_{\kappa,\mu}(z)$.
Thus the decaying solutions of the Schr\"odinger equation at $x = +\infty$ are the multiples of $e^{-x/2}W_{\kappa,\sqrt{\gamma-E}}(4\pi e^{x})$.  So the shooting function for $\kappa = \frac94$ is 
\begin{equation}
P(E) = 2 \pi^{-\frac14} W_{\frac94,\sqrt{\gamma-E}}(4\pi),
\end{equation}
where the prefactor is chosen to agree with the negative $E$ asymptotics of $\Xi$.
For the one-sided case we can now determine $\gamma = 0$, because using the integral representation 
\begin{equation}
W_{\kappa,\mu}(z) = \frac{z^{\mu+\frac12} 2^{-2\mu}}{\Gamma(\mu+\frac12-\kappa)} \int_1^\infty e^{-zt/2} (t-1)^{\mu-\frac12-\kappa} (t+1)^{\mu-\frac12+\kappa}\ dt
\end{equation}
one can obtain the asymptotic expansion
\begin{equation}
\partial_E \log P(E) \sim -\frac{1}{2\sqrt{-E}} \log \frac{\sqrt{-E}}{\pi} + \frac{7}{8E} + \frac{\gamma}{4(-E)^\frac32} \log\frac{\sqrt{-E}}{\pi e} + O((-E)^{-\frac32}),
\end{equation}
which is inconsistent with that for $\Xi$ unless $\gamma = 0$.

\begin{figure}[htbp] %  figure placement: here, top, bottom, or page
   \centering
   \includegraphics[width=5in]{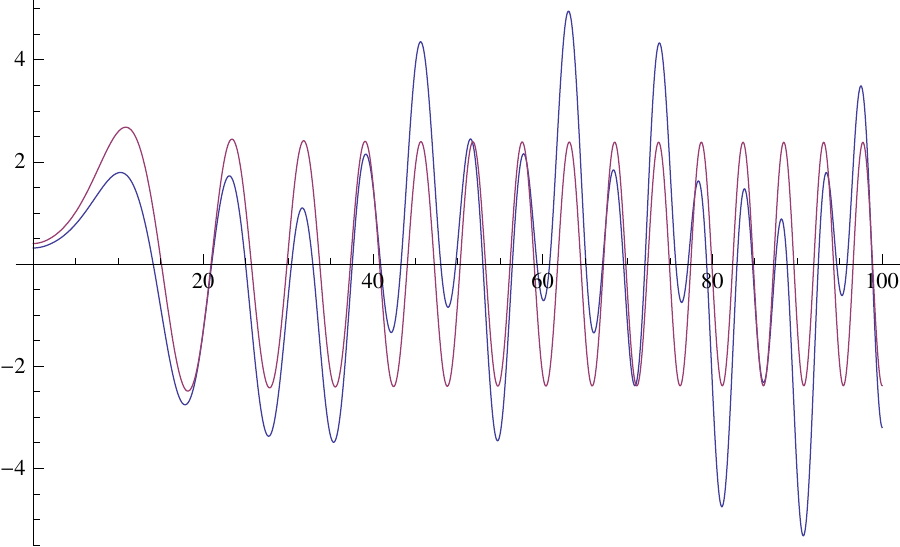} 
   \caption{Comparison of $\xi(\omega)$ (in blue) and $W_{\frac94,i \frac{\omega}{2}}(4\pi)$ (in red) as functions of $\omega$, scaled by $S$.}
   \label{fig:xiW}
\end{figure}

%\footnote{Whittaker's equation and functions are even in the order $\mu$ so the square root induces no singularity.}
Thus taking $\gamma=0$, Figure~\ref{fig:xiW} shows the match for positive $E = \omega^2/4$.
The fit is much better than for $K_{i\omega/2}(2\pi)$, but the oscillations are too regular.

\begin{figure}[htbp] %  figure placement: here, top, bottom, or page
   \centering
   \includegraphics[width=5in]{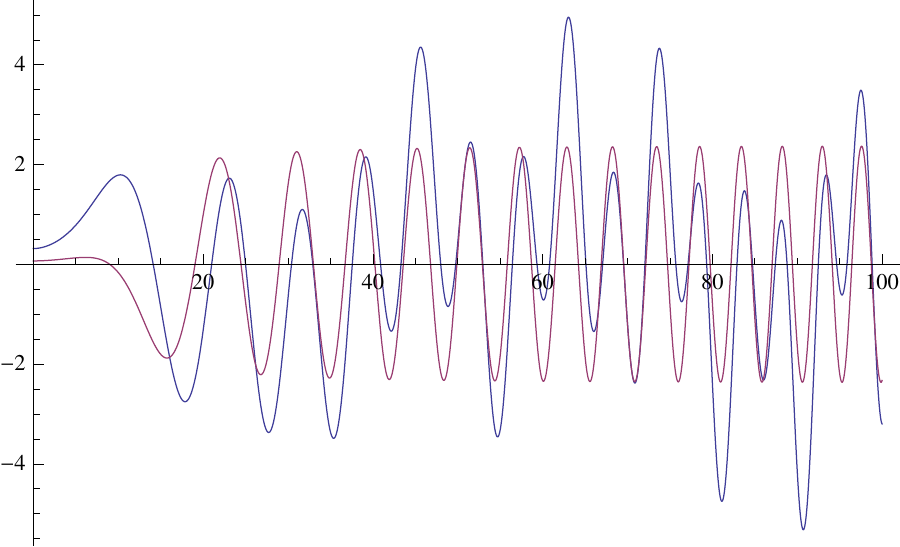} 
   \caption{Comparison of $\xi$ in blue with Polya's fake $\xi$-function $\xi^*$ in red (both scaled by $S$).}
   \label{fig:Polyafake}
\end{figure}

The fit is of similar quality to that of Polya's fake $\xi$-function $\xi^*$, which can be expressed \cite{Shi} as $4\pi^2 (K_{i\frac{\omega}{2}+\frac94}(2\pi) + K_{i\frac{\omega}{2}-\frac94}(2\pi))$, shown in Figure~\ref{fig:Polyafake}.
Polya proved all its zeroes are real by a non-spectral method.  The Whittaker function fits at least as well, however, and better for small $\omega$.  The addition of the correction $-6\pi (K_{i\frac{\omega}{2}+\frac54}(2\pi) + K_{i\frac{\omega}{2}-\frac54}(2\pi))$ to make Polya's second approximation, provides almost no visible improvement (so it is not shown).

Regularity of the oscillations of the shooting function is an unavoidable feature for smooth enough 1D Schr\"odinger operators.  Indeed, using Jeffreys' contribution to the LGWKB method \cite{Je}, namely a connection formula across the ``turning points'' where $V(x) = E$, valid if $V'' \ll |V'|^{4/3}$ for large $E$, the shooting function is the real part of a function $\tilde{P}$ with asymptotics in the one-sided case 
%[CHECK]
\begin{equation}
\tilde{P}(E) \sim 2(E-V(0))^{-\frac14} e^{-\int_0^a \sqrt{V(x)}\ dx} e^{\int_a^\infty \sqrt{V(x)-E}-\sqrt{V(x)}\ dx} e^{i(\int_0^a \sqrt{E-V(x)} \ dx - \frac{\pi}{4})},
\end{equation}
as $E \to \infty$, where $a(E)$ is defined by $V(a)=E$, assumed unique here.  The formula assumes that $V$ is bounded below and moreover that $V \ge 0$ (if not, one has to choose reference energy below or equal to the minimum).  

In the two-sided case (again assuming $V\ge 0$),
\begin{equation}
\tilde{P}(E) \sim -4 e^{-\int_{a_-}^{a_+} \sqrt{V(x)}\ dx} e^{\int_{I \setminus (a_-,a_+)} \sqrt{V(x)-E}-\sqrt{V(x)}\ dx} e^{i(\int_{a_-}^{a_+} \sqrt{E-V(x)}\ dx - \frac{\pi}{2})},
\end{equation}
where $a_\pm(E)$ are the left and right turning points, assumed unique.
Note that this result can be expressed purely in terms of the width function for the potential:
\begin{equation}
\tilde{P}(E) \sim -4 e^{-\int_0^E \sqrt{v} \ dw(v)} e^{\int_E^\infty \sqrt{v-E}- \sqrt{v} \ dw(v)} e^{i (\int_0^E \sqrt{E-v}\ dw(v) - \frac{\pi}{2})}.
\end{equation}

So to make a potential with irregular oscillations in its shooting function one would have to violate the condition $V'' \ll |V'|^{4/3}$.

\section{2D extensions}
\begin{figure}[htbp] %  figure placement: here, top, bottom, or page
   \centering
   \includegraphics[width=4in]{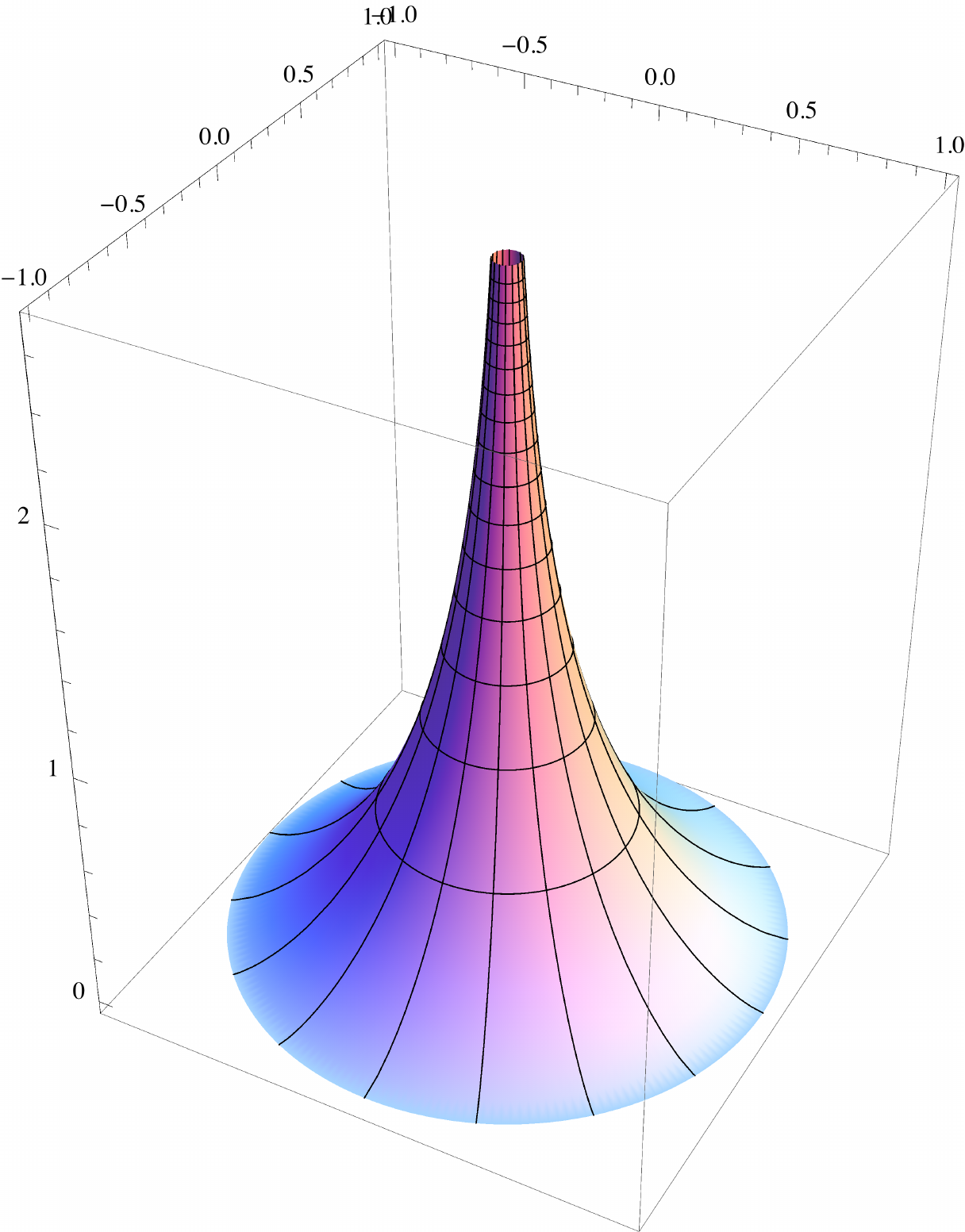} 
   \caption{An isometric embedding of a cusp of a surface of curvature $-1$ into Euclidean 3-space.}
   \label{fig:cusp}
\end{figure}

The appearance of the Whittaker function is suggestive of the Laplacian on a surface of curvature $-1$ with magnetic field $\frac94$, a ``cusp''  and a Dirichlet condition at a point.  A cusp in this sense is shown in Figure~\ref{fig:cusp}.\footnote{It is the image of the part of the upper half plane $z = x+iy$ with $y\ge 1/2\pi$, modulo $x \mapsto x+1$, under $(\frac{1}{2\pi y}\sin 2\pi x, \frac{1}{2\pi y}\cos 2\pi x, \log(2\pi y+\sqrt{4\pi^2 y^2-1})-\sqrt{1-(2\pi y)^{-2}})$.}
The Dirichlet point is taken at the point $z=i$, which is at height $\log(2\pi+\sqrt{4\pi^2-1})-\sqrt{1-(2\pi)^{-2}}$ $\approx 1.53738$ in the figure.

Indeed, for the magnetic Laplacian\footnote{I use the sign convention which makes it a positive operator.} $\Delta_{\frac94}$ on the upper half plane $z=x+iy, y>0$,  taking the gauge with magnetic 1-form $A = -\frac{9}{4y} dx$,
\begin{equation}
\Delta_{\frac94} = -y^2 (\partial_y^2 + \partial_x^2) + i\frac92 y \partial_x + \frac{81}{16},
\end{equation}
the solutions of $(\Delta_{\frac94}-\frac{85}{16}) \psi = E \psi$ at energy $E = \omega^2/4$, with periodicity under $z \mapsto z+1$ and decaying as $\Im z \to +\infty$ are
$$\psi(z) = \sum_{n>0} a_n W_{\frac94,i\frac{\omega}{2}}(4\pi n y) e^{-2\pi i n x} + \sum_{n<0} a_n W_{-\frac94,i\frac{\omega}{2}}(4\pi |n| y) e^{2\pi i |n| x}$$
for arbitrary coefficients $a_n$. 
If the surface has a Dirichlet condition at $z=i$ then we obtain the condition $P(E)=0$ with
\begin{equation}
P(E) = \sum_{n>0} a_n W_{\frac94,i\frac{\omega}{2}}(4\pi n) + \sum_{n<0} a_n W_{-\frac94,i\frac{\omega}{2}}(4\pi |n|),
\label{eq:sumW}
\end{equation}
of which the term $n=1$ is the 1D approximation of the previous section.
An attraction of this view is that the operator is truly complex rather than real, fitting with the numerically observed statistics of spacings between the Riemann zeroes.

One way this Dirichlet condition could arise, without actually putting in a nail at $z=i$, is if the surface is also invariant under $z \mapsto -1/z$ (producing the modular surface) and there is a flux string of strength $\pi$ at its fixed point $z=i$ so that 
\begin{equation}
\psi(-1/z) = - (-\frac{z}{\bar{z}})^\frac94 \psi(z).
\label{eq:fluxcond}
\end{equation}  
Taking $z \to i$ this would imply that $\psi(i)=0$.
Suppose there are other conditions which force the $(a_n)$ to be in a certain direction (possibly $\omega$-dependent).  Then the characteristic function for $\Delta_{\frac94}-\frac{85}{16}$ would be (\ref{eq:sumW}).

Now the Whittaker functions have asymptotic expression for large $\omega$ (e.g.~use \cite{Kaz})
\begin{equation}
W_{\kappa,i\frac{\omega}{2}}(Y) \sim e^{-\frac{\pi\omega}{4}} (\frac{\omega}{2})^{\kappa-\frac12} \sqrt{2Y} \cos(\frac{\omega}{2}\log \frac{2\omega}{Ye}+(\kappa-\frac12)\frac{\pi}{2}),
\label{eq:Wasymp}
\end{equation}
so 
\begin{eqnarray}
P(E) &\sim& \sum_{n>0} a_n e^{-\frac{\pi \omega}{4}} (\frac{\omega}{2})^\frac74 \sqrt{8\pi n} \cos(\frac{\omega}{2}\log \frac{\omega}{2\pi n e} + \frac78 \pi) \\
&+& \sum_{n<0} a_n e^{-\frac{\pi\omega}{4}} (\frac{\omega}{2})^{-\frac{11}{4}} \sqrt{8\pi |n|} \cos (\frac{\omega}{2}\log \frac{\omega}{2\pi |n| e} - \frac{11}{8}\pi).\nonumber
\end{eqnarray}

We compare this with the asymptotic expansion for $\xi$ discovered in Riemann's notes by Siegel.  Taking out the prefactor $f(\omega)$ it can be written as
\begin{equation}
- Z(\omega) = 2 \sum_{1\le N < \sqrt{\omega/2\pi}} N^{-\frac12} \cos \left(\arg \Gamma(\frac14+i\frac{\omega}{2}) - \frac{\omega}{2} \log \pi N^2 + \pi\right) + O(\omega^{-1/4}),
\end{equation}
and the argument of the cosine is asymptotically $\frac{\omega}{2}\log \frac{\omega}{2\pi N^2 e} + \frac78\pi$.
Actually, the expansion of $-Z$ is easy to derive by Laplace's saddle point method applied to the Fourier transform definition (\ref{eq:xiPhi}) of $\xi$, by deforming the contour of integration for each term of Polya's $\Phi$ to go through its saddle.  The real significance of Riemann's work was an explicit formula for the remainder term.

\begin{figure}[htbp] %  figure placement: here, top, bottom, or page
   \centering
   \includegraphics[width=5in]{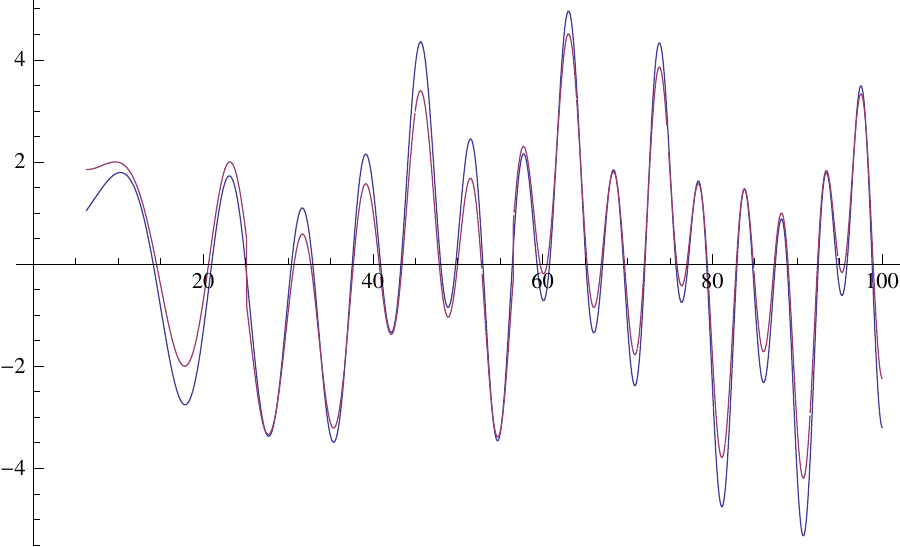} 
   \caption{$S\xi$ (in blue) and the Riemann-Siegel approximation for $-Z$ (in red) as functions of $\omega$, starting at $2\pi$.}
   \label{fig:RS}
\end{figure}

To see how good is the Riemann-Siegel asymptotic formula, Figure~\ref{fig:RS} shows the Riemann-Siegel approximation for $-Z$ and the function $S\xi$.
The Riemann-Siegel approximation jumps whenever $\omega$ passes through $2\pi N^2$ for $N$ integer, but the jumps are barely visible.  The fit is remarkably good, especially for the zeroes, though the amplitudes of the peaks and troughs are not accurate for this range of $\omega$.

The comparison between the Riemann-Siegel formula and the above $P(E)$ suggests that when $n$ is a square $N^2$ we choose $a_n \sim N^{-3/2} = n^{-3/4}$ (as $\omega \to \infty$), for other positive $n$ we choose $a_n$ anything substantially less than $n^{-3/4}$, and for $n<0$ we choose $a_n$ anything substantially less than $\omega^\frac92 |n|^{-\frac34}$.  There is no need to truncate to $N < \sqrt{\omega/2\pi}$ because the Whittaker functions $W_{\kappa,i\omega/2}(4\pi n)$ are exponentially negligible for $\omega \ll 4\pi n$.

\begin{figure}[htbp] %  figure placement: here, top, bottom, or page
   \centering
   \includegraphics[width=5in]{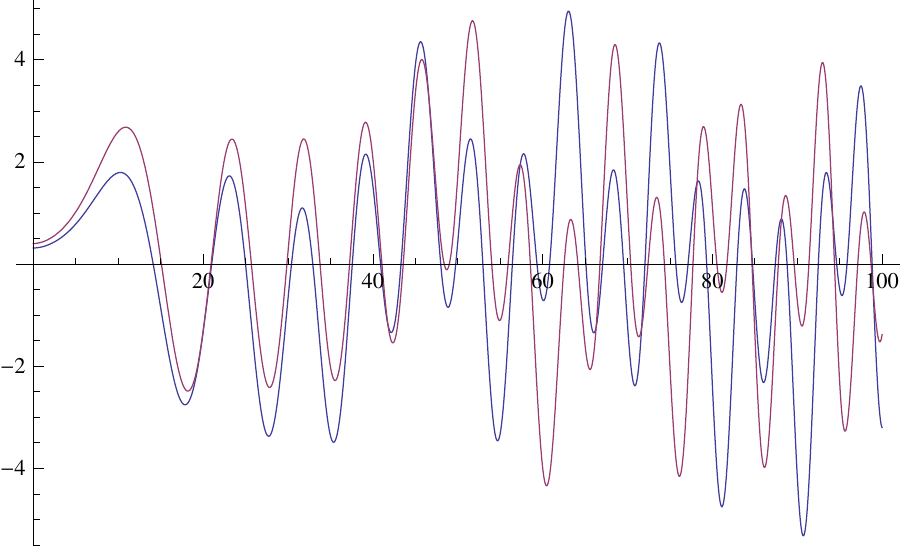} 
   \caption{Comparison of $\xi$ (in blue) with a sum of Whittaker functions (in red) (both scaled by $S$).}
   \label{fig:sumW}
\end{figure}

Taking $a_n = N^{-3/2}$ for $n = N^2$, and $0$ for the rest (truncating after the third term because the remaining terms are smaller than the resolution of the picture) produces Figure~\ref{fig:sumW}.  Unfortunately the fit is not good.
%The Whittaker function at $4\pi n^2$ has the asymptotic expression
%\begin{equation}
%2\pi^{-1/4} S(\omega) W_{\frac94, i\omega/2}(4\pi n^2) \sim n \cos (\arg \Gamma(\frac14+i\frac{\omega}{2}) - \frac{\omega}{2} \log \pi n^2)
%\label{eq:Wn}
%\end{equation}
%so it is natural to ask whether an equally good fit to $-Z$ would be obtained by
%$$2 \pi^{-1/4} S(\omega) \sum_{n=1}^{\sqrt{\omega/2\pi}} n^{-3/2} W_{\frac94,i\omega/2}(4\pi n^2)$$
%Figure~\ref{fig:sumW} shows the answer.  Unfortunately it does not fit well.
One problem is that the asymptotic expression (\ref{eq:Wasymp}) is attained non-uniformly in $Y$.  
Another is the negligibility of $W_{\kappa,i\omega/2}(4\pi n)$ for $\omega \ll 4\pi n$, which implies that the first three oscillations of $W_{9/4,i\omega/2}(4\pi)$ remain virtually unchanged on addition of the terms with $n = N^2$ for $N\ge 2$.
Perhaps the fit could be improved by adding subdominant terms for non-square $n$, but the choice would need to be based on some principle and I did not find a good one yet.  
A suggestion made by Yi-Chiuan Chen is to seek to determine the coefficients of $\xi$ in an expansion in Whittaker functions of the form (\ref{eq:sumW}) by finding a suitable inner product with respect to $\omega$.  Then one might understand how they could be enforced.
One could also make the coefficients $a_n$ depend on $\omega$, but the choice is limited by the requirement to keep $P(E)$ entire of order $\frac12$, so it seems to me the only freedom is to add $(A_n+B_n \omega^2 + C_n \omega^4) W_{-\frac94,i\frac{\omega}{2}}(4\pi |n|)$ for $n < 0$.

This is as far as I have reached, but I wish to point out some further difficulties.

Firstly, what sort of conditions could impose a direction on the coefficients $a_n$?  Perhaps one could take inspiration from \cite{Z} in which it is shown that for the non-magnetic Laplacian on the modular surface with condition that the integral along each of a certain set of closed geodesics is zero, the spectrum contains the Riemann zeroes.  This astonishing result would prove RH if one could show that with these conditions the operator is Hermitian, but it seems this has been tried and found hard.  Some analogously ``dense'' conditions would be required for the magnetic Laplacian to reduce from Weyl's law\footnote{the modular surface has area $\pi/3$; there are some corrections for cusps, e.g.~\cite{Mu}, and the only reference I found for cases with constant magnetic field assumes the cusps are plugged \cite{MT}, but I expect Weyl's law holds for our case too.} $N(E) \sim AE/4\pi$ for the number of eigenvalues less than $E$ for a surface of area $A$ 
to the required $\frac{\sqrt{E}}{\pi} \log \frac{\sqrt{E}}{\pi e}$.  But the magnetic geodesics are energy-dependent so they are unlikely to play a role.

Secondly, if we impose the Dirichlet condition at $z=i$ by a flux string on the modular surface we have a difficulty with Dirac quantisation of magnetic fields (which \cite{AKPS} attribute to Petersson), which says that the integral of the magnetic field over a surface must be a multiple of $2\pi$.  The modular surface has area $\pi/3$, so taking the flux string into account, the integral of the field is $\frac94 \frac{\pi}{3} - \pi = -\frac{\pi}{4}$.  This tells us that another flux string is required.  There could be one through the order-3 point of the modular surface, but its flux should be a multiple of $2\pi/3$ so that does not help us.  Or there could be one through the cusp, but that would change the periodicity condition on the wave function\footnote{the necessity of this can be seen from (\ref{eq:fluxcond}) evaluated at $z = (i\sqrt{3} \pm 1)/2$.} to $\psi(z+1) = e^{-i\pi/4} \psi(z)$,
and hence the Whittaker expansion should have $n+\frac18$ in place of $n$, which would not fit the Riemann-Siegel expansion.  Perhaps instead of a closed surface one should take one with boundary (a magnetic billiard domain).

\section{Comments}

Berry and Keating have been champions of the spectral approach to RH, e.g.~\cite{BK}.

Comtet and co-workers (e.g.~\cite{CGO}) and Avron and co-workers (e.g.~\cite{AKPS}) have been studying magnetic Laplacians on surfaces of curvature $-1$ for a long time, as a playground for complex Hermitian semiclassical mechanics and for mesoscopic quantum physics.  

I am not aware, however, that anyone has seen before that the good value to take for the magnetic field for RH is $\frac94$.  This assumes that it looks like a one-sided problem when approximated in 1D.
If the relevant 1D problem were two-sided (not necessarily symmetric) then, as I have pointed out, the relevant magnetic field would be $\frac74$.

\section*{Acknowledgements}
I am grateful to the University of Warwick for study leave for the academic year 2015/6, which allowed me time to develop these ideas.
The figures were produced using Mathematica.

\end{document}